\documentclass[12pt]{article}

\usepackage[margin=1in]{geometry}
\usepackage{amsmath,amssymb}
\usepackage{graphicx}
\usepackage{booktabs}
\usepackage{pgfplots}
\usepackage{hyperref}
\usepackage[round]{natbib}
\usepackage{adjustbox}
\usepackage{longtable}
\usepackage{booktabs}

\usepackage[most]{tcolorbox}
\usepackage[table]{xcolor}
\usepackage{tabularx}
\usepackage{array}

\usepackage{microtype}
\newcolumntype{Y}{>{\centering\arraybackslash}X}

\usepackage[most]{tcolorbox}

\newtcolorbox{activitybox}[1]{
    breakable,
    colback=gray!5,
    colframe=black!45,
    title=\textbf{Activity Handout: #1},
    fonttitle=\bfseries,
    boxrule=0.6pt,
    arc=2pt,
    left=8pt,
    right=8pt,
    top=6pt,
    bottom=6pt
}

\pgfplotsset{compat=1.18}

\title{Fast-Track ODE: A Flipped, Game-Based Summer Differential Equations Course}
\author{Chamila Malagoda Gamage\\
Department of Mathematics, University of Florida\\
\href{mailto:cgamage@ufl.edu}{cgamage@ufl.edu}}
\date{}

\begin{document}

\maketitle

\begin{abstract}
This article describes the implementation of a flipped, interactive, and game-based instructional design in an Elementary Differential Equations course taught during a six-week summer session. The course used real-time feedback, technology-supported visualization, collaborative problem solving, and mathematical games to support engagement, participation, and understanding in a fast-paced learning environment.
\end{abstract}


\section{Introduction}

Elementary Differential Equations is a course in which students are expected to connect calculus, algebra, modeling, and interpretation. Students learn several solution methods, but they are also expected to understand what differential equations represent and how they can be used to describe changing quantities. This can be difficult even in a regular semester. It becomes more challenging in a summer session, where the course is compressed into a much shorter period of time.

In a six-week summer course, attendance and engagement are especially important. Missing one class may mean missing an entire method or application. Students also need to keep attention for longer class meetings and move quickly from one topic to the next. Because of this, the course cannot depend only on traditional lecture. Students need frequent opportunities to practice, ask questions, check their understanding, and see why the material is useful.

In this paper, we describe a classroom-based redesign of a six-week Elementary Differential Equations course using a flipped, interactive, and game-based approach. Several in-class activities were introduced to make class time more active, provide frequent feedback, and help students connect differential equations with meaningful applications. We also discuss student feedback from the course, which suggested that students found the activities helpful for engagement, participation, exam preparation, and understanding.

\section{Background and Related Literature}

\subsection{Teaching Differential Equations in a Compressed Summer Course}

Teaching differential equations in a summer session creates a special instructional challenge. The course contains several topics that students often find new or abstract, including first-order equations, modeling, second-order equations, Laplace transforms, and systems. In a regular semester, students have more time to practice these ideas and revisit earlier material. In a six-week course, the pace is much faster, so students need to stay engaged from the beginning.

\citet{Luo2024} describes an intensive five-week summer differential equations course and emphasizes the importance of engagement, flexibility, and real-world applications. The present paper is related to this work because it also focuses on teaching differential equations in a compressed summer format. However, the course described here was an in-person six-week course built around flipped instruction, real-time feedback, applications, and game-based activities.

\subsection{Active Learning, Student-Centered Instruction, and Real-Time Feedback}

This course design is connected to active learning and student-centered instruction. Rather than using class time only for extended lecture, the course included frequent opportunities for students to solve problems, discuss ideas, answer questions, work with classmates, and explain their reasoning. This approach is consistent with active learning, where students are asked to participate directly in the learning process rather than only listen to the instructor \citep{Prince2004}. \cite{Freeman2014} also found that active learning in undergraduate STEM courses is associated with improved student performance compared with traditional lecture.

The design is also connected to formative assessment. In a formative view of assessment, student responses are used not only to evaluate learning after instruction, but also to guide teaching while learning is still taking place \citep{BlackWiliam1998}. Classroom assessment techniques can give instructors low-stakes information about student understanding while there is still time to respond \citep{AngeloCross1993}. This was especially important in a six-week course, where misunderstandings needed to be identified, and act quickly.

\cite{ClineLomen2010} give a differential-equations example of this idea through classroom voting. Their approach uses multiple-choice questions, student discussion, voting, and explanation to involve students more actively in the lesson. The online quizzes, polls, and quick questions used in the present course served a similar purpose: they made student thinking visible, encouraged participation, and helped the instructor decide when to move forward or revisit an idea.

\subsection{Modeling, Visualization, and Applications}

Differential equations are naturally connected to modeling, but students may still experience the course mainly as a collection of solution methods. They may learn how to solve an equation without fully understanding what the equation represents or why the solution matters. For this reason, modeling, visualization, and applications were important parts of the course design.

\cite{Huber2010} describes the use of interdisciplinary application projects in an ordinary differential equations course. This work is relevant because it emphasizes real-world problems, applications from other disciplines, and the use of modeling to make solution methods more meaningful. The present course used the same general idea, but adapted it to a fast-paced summer setting by building applications directly into in-class activities.

Technology also supported this goal. Desmos was used to help students visualize solution curves, parameters, and graphical behavior. This helped students see that differential equations are not only symbolic objects; they describe motion, growth, decay, cooling, oscillation, and other changing systems.

\subsection{Game-Based Activities as Structured Mathematical Practice}

The game-based activities in this course were designed to support mathematical practice, not to replace it. Game-based learning can create learning environments in which students engage with content through goals, feedback, interaction, and problem solving \citep{PlassHomerKinzer2015}. In this course, each activity had a specific mathematical purpose, such as reviewing for exams, interpreting a model, practicing Laplace transforms, or explaining a solution method. The game structure made the work more engaging, but the focus remained on the mathematics.

These activities are also connected to active learning. Students were not only watching the instructor solve problems; they were solving, discussing, explaining, checking, and defending their reasoning. This is important in differential equations because students need to choose methods, interpret results, and connect answers back to the original problem. Prior research on serious games suggests that well-designed games can support learning when they are connected to clear instructional goals \citep{Wouters2013}. In the present course, the games were used as structured opportunities for practice, feedback, and collaboration.

Overall, this paper builds on prior work on summer differential equations instruction, classroom feedback, modeling, visualization, active learning, and game-based learning. Its main contribution is to describe how these ideas can be combined into a practical classroom structure for a fast-paced six-week differential equations course.

\section{Course Context}

The course was an Elementary Differential Equations course taught during a six-week summer session at a large public university in the United States. The class had 30 students representing 17 majors, including economics, engineering, psychology, music, and other fields. This diversity of academic backgrounds was important because students came to the course with different mathematical preparation, interests, and reasons for taking differential equations.

The summer course covered the same main topics as the regular semester version of the course. These included first-order differential equations, applications and modeling, second-order differential equations, and Laplace transforms. The class met five days per week for 75 minutes each day. Students also completed online homework, in-person quizzes, and exams as part of the regular course structure.

By nature, the course was fast-paced. A new topic was introduced almost every day, leaving limited time for students to fall behind and catch up later. For this reason, regular attendance, active participation, and frequent checks of understanding were especially important. 

\section{Instructional Design and Activities}

The instructional design was built around frequent interaction and student participation. On the first day of class, students responded to a word-cloud prompt asking what came to mind when they heard the phrase ``differential equations.'' This served as an icebreaker and opened a class discussion about what differential equations are, where they appear, and why they are useful. Throughout the course, when a new idea, method, or concept was introduced, students often responded to a quick poll, quiz, or short question using tools such as Mentimeter. These real-time responses allowed students to participate immediately and gave the instructor a quick view of student understanding. When the responses showed confusion, the instructor could slow down, revisit an example, or adjust the explanation before moving forward. Before each exam, one or two class meetings were used for review, and these review sessions were often organized as games so that students could review important topics while staying active and engaged. In addition to these regular interactive elements, several other activities were designed for specific course topics.

\subsection{Real-Time Polls, Quizzes, and Word Clouds}

\begin{figure}[htbp]
    \centering
     \begin{adjustbox}{frame, margin=6pt}
     \includegraphics[width=1\textwidth]{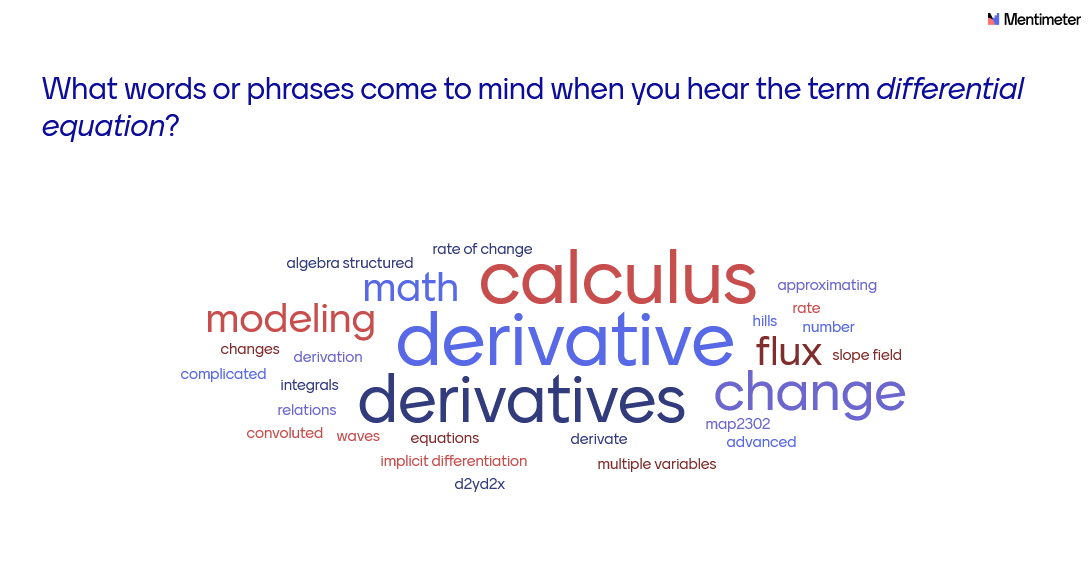}
    \end{adjustbox}
    \caption{Example of a word cloud used during class.}
    \label{fig:mentimeter-examples}
\end{figure}

\begin{figure}[htbp]
    \centering
     \begin{adjustbox}{frame, margin=6pt}
     \includegraphics[width=1\textwidth]{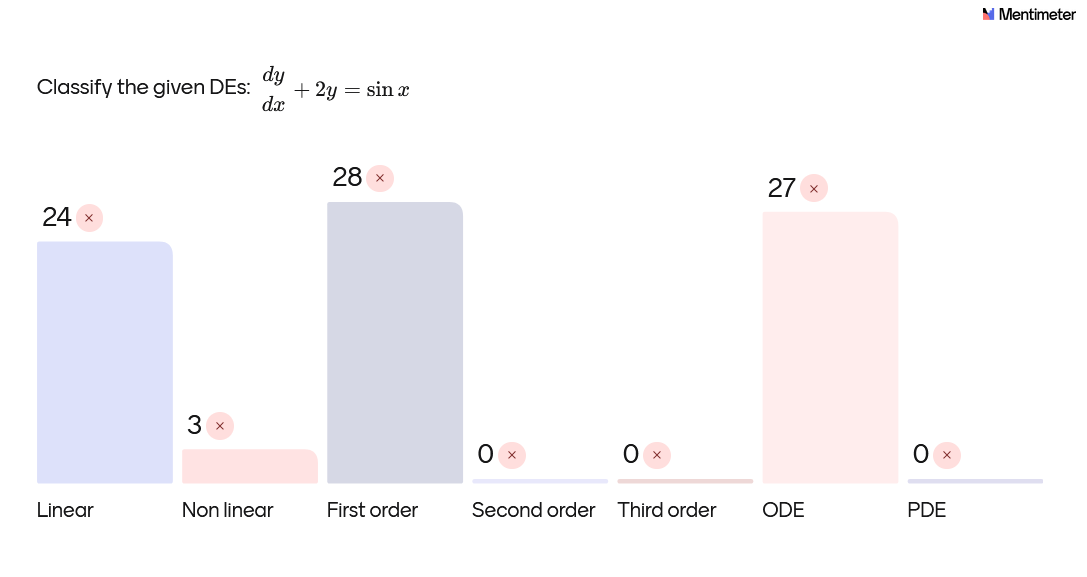}
    \end{adjustbox}
    \caption{Example of a real-time quiz used during class.}
    \label{fig:mentimeter-examples2}
\end{figure}

Real-time polls, quizzes, and word clouds were used throughout the course to make student thinking visible during class. Figures~\ref{fig:mentimeter-examples} and \ref{fig:mentimeter-examples2} show examples of a word cloud and a quiz question given through Mentimeter. Students could quickly scan a QR code, submit their answers, and see the class responses displayed immediately. This allowed the instructor to continue the class discussion based on students' actual responses and to adjust the explanation when many students showed confusion.

\newpage

\subsection{Connecting Differential Equations to Students' Fields}

The class included students from 17 majors, and students had different reasons for taking differential equations. Because of this, the applications part of the course was used to emphasize that differential equations are a versatile tool with connections to many fields. When introducing mathematical modeling, the instructor showed examples from several disciplines represented in the class. Many students were surprised to see direct connections between differential equations and their own areas of study, which helped them view the course as more than a collection of formulas and solution methods.
\vspace{1cm}

\begin{table}[htbp]
    \centering
    \small
    \renewcommand{\arraystretch}{1.5}
    \begin{tabular}{p{0.28\textwidth} p{0.62\textwidth}}
        \toprule
        \textbf{Field} & \textbf{Example Application} \\
        \midrule
        Economics & Capital growth and depreciation in the Solow--Swan growth model. \\
        Mechanical Engineering & Motion of a mass--spring--damper system. \\
        Material Science & Diffusion of concentration through a material. \\
        Mathematics & Heat flow in a rod with boundary and initial conditions. \\
        Nuclear Engineering & Neutron density in point reactor kinetics. \\
        Health Science & Spread and recovery in an epidemic model. \\
        Astronomy & Motion of two bodies under gravity. \\
        Aerospace Engineering & Rocket velocity as mass changes during flight. \\
        Computer Engineering & Voltage across a capacitor in an RC circuit. \\
        Physics & Evolution of a quantum wavefunction. \\
        Environmental Engineering & Transport and removal of a pollutant. \\
        Biomedical Engineering & Drug concentration in the body. \\
        Chemical Engineering & Reactant concentration in a continuously stirred tank reactor. \\
        Statistics & Mean-reverting random behavior in an Ornstein--Uhlenbeck process. \\
        Music & Vibration of a string and wave motion. \\
        Philosophy & Replicator dynamics for competing strategies or ideas. \\
        \bottomrule
    \end{tabular}
    \caption{Examples used to show how differential equations appear across different fields.}
    \label{tab:de-fields}
\end{table}

\newpage

\subsection{Courtroom Activity: Newton's Law of Cooling}

Several applications of first-order differential equations were included in the course. This was important because, without applications, students may experience the topic as memorizing formulas and completing long calculations. The goal was to help students see that a differential equation can model a real situation and that the solution has meaning in context. When Newton's Law of Cooling was introduced, students worked in groups on a courtroom-style activity. The activity was designed as a trial, where students used cooling-law calculations as evidence. Students took the roles seriously, worked through the mathematics carefully, and enjoyed presenting their arguments. Although the activity had a playful structure, the main purpose was to practice setting up and solving a first-order model and interpreting the result.

\vspace{1cm}

\noindent
\begin{activitybox}{Murder at the Mansion: A Cooling-Law Defense}

\medskip

\noindent
\textbf{Scenario Description.}
Late one evening, the police are summoned to the mansion of Mr.\ Green and discover him slumped in his private study. A crime scene investigator records his body temperature twice, once at 20:00 and again at 21:00, while the room remains at a steady \(20^\circ\mathrm{C}\). Four suspects were seen near the study and have provided their alibis.

\begin{itemize}
    \item \textbf{Anna} claims she was in the study from 18:30 to 18:50.
    \item \textbf{Ben} says he stayed from 18:45 to 19:15.
    \item \textbf{Emma} reports being there briefly between 19:00 and 19:08.
    \item \textbf{Dan} insists he was at the gym until 20:00, but his fitness tracker shows he left at 19:00 and lingered nearby until 19:30.
\end{itemize}

\noindent
Now the suspects have been called to court. Each suspect's legal team will use the cooling-law calculations to defend their client, while the jury will verify the mathematics, ask questions, and deliver a verdict.

\medskip

\noindent
\textbf{Newton's Law of Cooling}
\[
    \frac{dT}{dt} = k(T_{\mathrm{env}}-T),
    \qquad
    T(t) = T_{\mathrm{env}} + (T_0-T_{\mathrm{env}})e^{-kt}.
\]

\noindent
\textbf{Given Data:}
\begin{itemize}
    \item \(T_0 = 37^\circ\mathrm{C}\) \quad normal body temperature,
    \item \(T_{\mathrm{env}} = 20^\circ\mathrm{C}\) \quad room temperature,
    \item \(T(20{:}00)=32^\circ\mathrm{C}\),
    \item \(T(21{:}00)=28^\circ\mathrm{C}\).
\end{itemize}

\noindent
\textbf{Group Roles:}
\begin{itemize}
    \item \textbf{Legal Teams} one team per suspect:
    \begin{itemize}
        \item \textit{Math Experts} perform the cooling-law calculations.
        \item \textit{Presenter} explains the results and timeline.
        \item \textit{Creative Strategist} builds arguments to defend the suspect.
    \end{itemize}

    \item \textbf{Jury Team:}
    \begin{itemize}
        \item \textit{Math Verifiers} verify the calculations.
        \item \textit{Questioner} asks clarifying questions from each defense team.
        \item \textit{Verdict Officer} tallies votes and announces the ``Best Defense.''
    \end{itemize}
\end{itemize}

\noindent
\textbf{Activity Steps:}
\begin{enumerate}
    \item \textbf{All teams:}
    \begin{enumerate}
        \item Find the cooling constant \(k\).
        \item Compute the time since death at 20:00.
        \item Determine the clock-time of death.
        \item Identify the most likely suspect.
    \end{enumerate}
    Students may use a calculator, but no additional sources or aids.

    \item \textbf{Defense Teams:}
    \begin{itemize}
        \item Show the calculation clearly.
        \item Craft two arguments for why your client could not have committed the crime at that time.
        \item Prepare a two-minute opening statement to defend your client.
    \end{itemize}

    \item \textbf{Jury Team:}
    \begin{itemize}
        \item Ask one clarifying question per team.
        \item Vote on the most convincing defense.
    \end{itemize}
\end{enumerate}
\end{activitybox}

\subsection{Jeopardy-Style Exam Review Games}

Reviewing for exams was another important part of the course design. In a fast-paced summer course, review days need to cover many topics while also keeping students attentive and engaged. To address this, some exam review sessions were organized as Jeopardy-style games. Students worked in groups and selected problems based on topic categories and point values. After a group selected a problem, students worked together to solve it, explained their reasoning, and then compared their solution with the answer. This format allowed the class to review a large amount of material while keeping the review session active and collaborative.

Students were enthusiastic during these review sessions. The game format encouraged students to attempt problems, discuss strategies with classmates, and recall methods from earlier lessons. By the end of the activity, students had reviewed the main exam topics in a structured way and many reported that they felt more prepared for the exam.

Table~\ref{tab:jeopardy-board} shows an example of a Jeopardy-style review board used before an exam. Figure~\ref{fig:jeopardy-example} shows an example of a question and answer pair from the activity.

\begin{table}[htbp]
    \centering
    \renewcommand{\arraystretch}{1.5}
    \setlength{\tabcolsep}{4pt}
    \begin{tabular}{ccccc}
        \toprule
        Introduction & Initial Value Problems & Slope Fields & Separable Equations & Integrating Factors \\
        (1.1) & (1.2) & (1.3) & (2.2) & (2.3) \\
        \midrule
        100 & 100 & 100 & 100 & 100 \\
        200 & 200 & 200 & 200 & 200 \\
        300 & 300 & 300 & 300 & 300 \\
        400 & 400 & 400 & 400 & 400 \\
        500 & 500 & 500 & 500 & 500 \\
        \bottomrule
    \end{tabular}
    \caption{Example Jeopardy-style review board.}
    \label{tab:jeopardy-board}
\end{table}

\begin{figure}[htbp]
    \centering
    \fbox{
    \begin{minipage}{0.85\textwidth}
        \textbf{Categorization and Modeling -- 100}

        \medskip

        Classify the differential equation
        \[
            y'' + 3xy' + 2y = \sin x
        \]
        as linear or nonlinear.

        \medskip

        \textbf{Answer:} Linear.
    \end{minipage}
    }
    \caption{Example question and answer from a Jeopardy-style review game.}
    \label{fig:jeopardy-example}
\end{figure}

A second review activity used the same question-board structure, but instead of earning points, students earned items for an ``Escape the Island'' challenge. This showed that review games do not always need to be score-based. The same mathematical questions can be organized around tools, clues, or a final mission, while still keeping the focus on practice, discussion, and exam preparation.

\begin{table}[htbp]
    \centering
    \small
    \renewcommand{\arraystretch}{1.25}
    \setlength{\tabcolsep}{3pt}
    \begin{tabularx}{\textwidth}{
        >{\raggedright\arraybackslash}p{0.23\textwidth}
        >{\raggedright\arraybackslash}X
        >{\raggedright\arraybackslash}X
        >{\raggedright\arraybackslash}X
        >{\raggedright\arraybackslash}X
    }
        \toprule
        \textbf{Topic} & \textbf{200: Comfort} & \textbf{300: Gear} & \textbf{400: Knowledge} & \textbf{500: Escape} \\
        \midrule
        Integrating Factors & Granola bar & Compass & Hand-drawn map & Flare gun \\
        Exact Equations & Pair of socks & Can opener & Survival book & Radio transmitter \\
        Substitution Methods & Umbrella & Duct tape & Tracking guide & Raft repair kit \\
        Modeling & Coin & Spoon & Fire-starting instructions & GPS chip \\
        Second-Order ODEs & Deck of cards & Glow stick & Morse code chart & Signal mirror \\
        \bottomrule
    \end{tabularx}
    \caption{Example tool board from the Escape the Island review activity.}
    \label{tab:escape-island-tools}
\end{table}
\newpage

\subsection{Laplace Black Box}

The Laplace transform can be difficult for students to comprehend at first, because it asks them to move between two different spaces. A differential equation in the \(t\)-domain is transformed into an algebraic equation in the \(s\)-domain, solved there, and then transformed back into a solution in the \(t\)-domain. To emphasize this idea, students completed a partner activity called the Laplace Black Box. The activity was meant to imitate the idea of a black box: students see the input and output, but they must understand the process well enough to move between the two representations. Students worked in pairs. One student transformed a differential equation and solved for \(Y(s)\), while the other student received only the final expression for \(Y(s)\) and used the inverse Laplace transform to recover \(y(t)\). The pair then worked together to verify that the solution satisfied the original differential equation and initial conditions.

\vspace{1cm}

\noindent
\begin{activitybox}{Laplace Black Box -- Partner Activity}

\medskip

\noindent
\textbf{Why ``Laplace Black Box''?}\\

\noindent
When we apply the Laplace transform to a differential equation, we treat the process like a black box. We encode the original equation into an algebraic expression in terms of \(Y(s)\), solve it more easily in the \(s\)-domain, and then decode it back into a solution \(y(t)\) in the time domain. In this activity, you and your partner will each explore one half of this process and then come together to verify that the solution works.

\medskip

\noindent
\textbf{Instructions}

\begin{enumerate}
    \item You will receive a differential equation with initial conditions.
    
    \item Take the Laplace transform of both sides of the equation and solve for \(Y(s)\). Simplify your answer as much as possible.
    
    \item Write only your final expression for \(Y(s)\) in the box provided.
    
    \item Exchange your \(Y(s)\) expression with your partner. You will receive a \(Y(s)\) expression from them as well.
    
    \item Take the inverse Laplace transform of your partner's \(Y(s)\) to recover the solution \(y(t)\).
    
    \item Work together to verify that the resulting \(y(t)\) satisfies the original differential equation and initial conditions.
\end{enumerate}

\medskip

\noindent
\textbf{Partner 1: Solve for \(Y(s)\)}

\medskip

\noindent
Write your differential equation and initial conditions here:

\[
    \rule{0.9\textwidth}{0.4pt}
\]

\noindent
Now take the Laplace transform of both sides and solve for \(Y(s)\):

\[
    Y(s) = \rule{0.6\textwidth}{0.4pt}
\]

\medskip

\noindent
\textbf{Pass This to Your Partner}

\medskip

\noindent
Write only your final expression for \(Y(s)\) below. This is what your partner will use to find \(y(t)\):

\[
    Y(s) = \rule{0.6\textwidth}{0.4pt}
\]

\medskip

\noindent
Use the expression for \(Y(s)\) above to find the corresponding solution \(y(t)\) using the inverse Laplace transform:

\[
    y(t) = \rule{0.6\textwidth}{0.4pt}
\]

\medskip

\noindent
Finally, verify with your partner that the solution satisfies the original differential equation and initial conditions.
\end{activitybox}

\subsection{Laplace Vault}

After students learned several Laplace transform formulas, they needed repeated practice using those formulas in different forms. The Laplace Vault activity was designed as a team-based puzzle to give students this practice in a more engaging way. Students first solved several inverse Laplace transform problems. After their answers were checked, they solved a short initial value problem using Laplace transforms. Finally, they used a decoder table to match their time-domain answers to letters and unscramble a final password. The purpose of the activity was to help students practice recognition of Laplace transform forms, inverse transforms, shifting, and IVP solving, while keeping the practice collaborative and active.

\vspace{1cm}

\noindent
\begin{activitybox}{Laplace Vault -- Final Challenge}

\medskip

\noindent
\textbf{Mission Briefing.}
You are stuck in the \(s\)-domain. To escape and return to the time-domain, your team must complete two tasks:

\begin{enumerate}
    \item Solve eight inverse Laplace transform problems. Each answer corresponds to a time-domain function.
    \item Once your instructor verifies that all eight answers are correct, you will receive a short IVP. If you solve it correctly, you will receive the decoder table.
\end{enumerate}

\noindent
Use the decoder table to match each of your time-domain expressions to a letter. Unscramble the letters to discover the final password and escape the vault.

\medskip

\noindent
\textbf{Phase 1: Inverse Laplace Transform Problems}

\begin{enumerate}
    \item
    \[
        \mathcal{L}^{-1}\left\{\frac{6}{(s+1)^2}\right\}
    \]

    \item
    \[
        \mathcal{L}^{-1}\left\{\frac{2}{s+2}\right\}
    \]

    \item
    \[
        \mathcal{L}^{-1}\left\{\frac{8s+16}{s(s+2)}\right\}
    \]

    \item
    \[
        \mathcal{L}^{-1}\left\{\frac{1}{s+1}\right\}
    \]

    \item
    \[
        \mathcal{L}^{-1}\left\{\frac{3}{s^2+1}\right\}
    \]

    \item
    \[
        \mathcal{L}^{-1}\left\{\frac{7}{s^2+4}\right\}
    \]

    \item
    \[
        \mathcal{L}^{-1}\left\{\frac{5}{s^2}\right\}
    \]

    \item
    \[
        \mathcal{L}^{-1}\left\{\frac{4e^{-3s}}{s+1}\right\}
    \]
\end{enumerate}

\medskip

\noindent
\textbf{Phase 2: Final Lock}

\medskip

\noindent
Solve the following initial value problem using Laplace transforms:
\[
    y''+y=e^{-t}, \qquad y(0)=0, \qquad y'(0)=0.
\]

\medskip

\noindent
\textbf{Phase 3: Decoder Table}

\begin{center}
\renewcommand{\arraystretch}{1.2}
\begin{tabular}{|c|l||c|l|}
\hline
\textbf{Letter} & \textbf{Time-Domain Expression}
&
\textbf{Letter} & \textbf{Time-Domain Expression}
\\
\hline
F & \(6te^{-t}\) & X & \(e^{-3t}\) \\
A & \(3\sin(t)\) & Y & \(3\sin(3t)\) \\
N & \(4e^{-(t-3)}u(t-3)\) & B & \(4e^{t-3}u(t-3)\) \\
T & \(e^{-t}\) & C & \(\cos(t)\) \\
M & \(8\) & K & \(e^{-t}\cos(t)\) \\
R & \(2e^{-2t}\) & O & \(\frac{7}{2}\sin(2t)\) \\
Q & \(7\sin(2t)\) & U & \(7\cos(2t)\) \\
V & \(8t\) & S & \(5t\) \\
H & \(te^{-2t}\) & L & \(t^2\) \\
Z & \(u(t-1)\sin(t-1)\) & D & \(\cos(2t)\) \\
E & \(6e^{-t}\) & & \\
\hline
\end{tabular}
\end{center}

\medskip

\noindent
\textbf{Your Letters.}
Match your answers to the decoder table above.

\begin{center}
\renewcommand{\arraystretch}{1.8}
\begin{tabular}{|c|c|c|c|c|c|c|c|}
\hline
1 & 2 & 3 & 4 & 5 & 6 & 7 & 8 \\
\hline
\hspace{2em} & \hspace{2em} & \hspace{2em} & \hspace{2em} &
\hspace{2em} & \hspace{2em} & \hspace{2em} & \hspace{2em} \\
\hline
\end{tabular}
\end{center}

\medskip

\noindent
\textbf{Unscrambled Secret Word:}
\rule{0.45\textwidth}{0.4pt}

\medskip

\noindent
\textit{Hint: Some letters may appear more than once in the secret word.}
    
\end{activitybox} 
\section{Student Feedback}

At the end of the course, students completed a short feedback survey about the in-class activities. The goal of the survey was not to make a formal statistical claim, but to understand how students experienced the activities and whether they found them useful. Fourteen students responded to the survey.

The overall response was very positive. In particular, students agreed that the activities helped them understand differential equations more effectively than traditional lectures, made them more interested in mathematics, motivated them to attend class regularly, helped them prepare for exams and quizzes, encouraged active participation, and should be included in future math classes.

The open-ended responses also showed that students found the activities useful and memorable. Several students connected the activities to applications of mathematics. For example, one student wrote that the courtroom trial ``helped me understand more clearly the practical applications of many mathematical formulas such as Newton's Law of Cooling.'' Another student described the courtroom activity as ``a ton of fun,'' while a different student wrote that ``the Jeopardy games were an amazing idea.'' One response summarized the value of the more creative activities by saying, ``the more out there the activity was the easier it was for me to engage with it.''

Overall, the student feedback suggests that the activities were well received. Students especially valued activities that helped them review for exams, work with classmates, and see how differential equations can be used in real situations. These responses support the main goal of the course design: to make a fast-paced summer course more active, engaging, and meaningful while keeping the focus on mathematical learning.

\section{Conclusion}

This paper described a classroom-based redesign of a six-week Elementary Differential Equations course. The main goal was to make a fast-paced summer course more active and meaningful, while still keeping the mathematical content at the center. The activities gave students regular opportunities to participate, practice, review, and see how differential equations can be used in different contexts.

The implementation also showed that this kind of course design requires some early planning. Instructors need time to prepare the activities, organize materials, form groups, and make sure that all students are included. This is especially important in a summer course, where the pace is fast and there is little extra time during the semester. Some students may also need encouragement to participate fully in group activities.

Overall, the redesign was successful. Student feedback suggested that the activities made the course more engaging and helped students connect with the material. For an intensive differential equations course, these activities provided a practical way to support attendance, participation, exam preparation, and student interest. More importantly, they helped students experience differential equations not only as a set of methods to memorize, but as a useful and active subject.

\bibliographystyle{apalike}
\bibliography{references}
\end{document}